\documentclass{article}%
\usepackage{amsmath,amsfonts,amssymb}

\usepackage{graphicx}
\usepackage{setspace}
\usepackage{hyperref}
\usepackage{mathtools}
\DeclarePairedDelimiter\floor{\lfloor}{\rfloor}
\usepackage{latexsym}
\usepackage{psfrag}
\usepackage{subfigure}
\usepackage{authblk}

\hypersetup{colorlinks=true, linkcolor=blue, citecolor=red, urlcolor=red}
\newtheorem{theorem}{Theorem}

\newtheorem{lemma}[theorem]{Lemma}

\newtheorem{proposition}[theorem]{Proposition}

\newenvironment{proof}[1][Proof]{\noindent\textbf{#1.} }{\ \rule{0.5em}{0.5em}}

\begin{document}
\doublespacing

\title{\textbf{A Bayesian Nonparametric Estimation to Entropy}}   


\author[1]{Luai Al-Labadi\thanks{{\em Corresponding author:} luai.allabadi@utoronto.ca}}

\author[1]{Vishakh Patel\thanks{vishakh.patel@mail.utoronto.ca}}

\author[1]{Kasra Vakiloroayaei\thanks{k.vakiloroayaei@mail.utoronto.ca}}

\author[1]{Clement Wan\thanks{cmclement.wan@mail.utoronto.ca}}

\affil[1]{Department of Mathematical and Computational Sciences, University of Toronto Mississauga, Mississauga, Ontario L5L 1C6, Canada.}

\date{}
\maketitle

\pagestyle {myheadings} \markboth {} {BNP Estimation to Entropy}

\begin{abstract}
A Bayesian nonparametric estimator to  entropy is proposed. The derivation  of the new estimator relies on using the Dirichlet process and adapting the well-known frequentist estimators of Vasicek (1976) and Ebrahimi, Pflughoeft and Soofi (1994).  Several theoretical properties, such as consistency,  of the proposed estimator are obtained. The quality of the proposed estimator has been investigated through several examples, in which it exhibits excellent performance.

\par
 \vspace{9pt} \noindent\textsc{Keywords:}  Bayesian non-parametric, Dirichlet process, Entropy, Kullback-Leibler divergence, Model checking.

 \vspace{9pt}

\end{abstract}

\pagebreak
\section{Introduction}

The concept of  (differential) \emph{entropy} was introduced in Shannon  (1948). Since then, entropy has been one of the most interesting areas with endless applications in many fields such as thermodynamics, communication theory,  computer science, biology, economic, mathematics and statistics (Cover and Thomas, 2006). The  entropy of a  continuous cumulative distribution function (cdf) $P$ with a probability density function $p$ (with respect to Lebesgue measure) is defined as
\begin{equation}
H(P)=-\int_{-\infty}^{\infty} p(x)\log p(x) dx=-E_P\left[\log p(x)\right]. \label{entropy}%
\end{equation}

From practical perspective, one must estimate (\ref{entropy}) from the data, which is not a trivial task. Various frequentist procedures for the estimation of entropy are offered in the literature. Among several estimators, due to its simplicity, Vasicek's  (1976) estimator has been the most common and the widely used one. Vasicek  (1976) noticed that (\ref{entropy}) can be written as
\begin{equation*}
H(P)=-\int_{0}^{1} \log \left(\frac{d}{dt}P^{-1}(t)\right) dt. \label{derivative}
\end{equation*}
Thus, $H(P)$ is estimated by using estimates of the derivative of inverse of the distribution function on the sample points. Specifically, if  $x=(x_{1},\ldots,x_{n})$ is a sample from a distribution $P$, then, at each sample point $x_i$, the derivative of $P^{-1}(t)$ is estimated by the slope defined by
\begin{equation}
\frac{x_{(i+m)}-x_{(i-m)}}{F_n(x_{(i+m)})-F_n(x_{(i-m)})}=\frac{x_{(i+m)}-x_{(i-m)}}{\frac{i+m}{n}-\frac{i-m}{n}}=\frac{x_{(i+m)}-x_{(i-m)}}{2m/n}, \label{slope}
\end{equation}
where $F_n$ is the empirical  distribution function. Consequently,  Vasicek  (1976) estimator is given by
\begin{align}
H^{V}_{m,n}=n^{-1}\sum_{i=1}^n  \log \left(\frac{x_{(i+m)}-x_{(i-m)}}{2m/n}\right), \label{Vasicek}
\end{align}
where  $m$, called the window size, is a positive integer smaller than $n/2$ and  $x_{(1)} \le x_{(2)} \le \cdots \le x_{(n)}$ are the order statistics of $x_{1}, x_{2}, \ldots, x_{n}$ with $x_{(i)}=x_{(1)}$ if $i<1$, $x_{(i)}=x_{(n)}$ if $i>n$. Vasicek  (1976) showed that $H^{V}_{m,n} \overset{p} \to H(P)$, where $\overset{p} \to$ denotes convergence in probability. Ebrahimi, Pflughoeft and Soofi (1994) noticed that (\ref{Vasicek}) does not give the correct formula for the slope when $i \leq m$ or $i \ge n-m+1$. They proposed the following modification to  (\ref{Vasicek}):
\begin{align}
H^{EPS}_{m,n}=n^{-1}\sum_{i=1}^n  \log \left(\frac{x_{(i+m)}-x_{(i-m)}}{c_im/n}\right), \label{EPS}
\end{align}
where
\begin{align}
c_{i} &= \left\{
        \begin{array}{ll}
            \frac{m+i-1}{m} & \quad  1\leq i \leq m \\
 \label{EPS1}           2 & \quad  m+1 \leq i \leq N-m\\
            \frac{N+m-i}{m} & \quad  N-m+1 \leq i \leq N\\
        \end{array}
    \right..
\end{align}
They also showed that $H^{EPS}_{m,N} \overset{p} \to H(P)$.

Other nonparametric frequentist  estimators of entropy  includes, among others, the work of  van Es (1992), Correa (1995),  Wieczorkowski and Grzegorzewski (1999),  Alizadeh Noughabi (2010), Alizadeh Noughabi  and Arghami  (2010), Bouzebda, Elhattab, Keziou and Lounis (2013) and Al-Omari (2014, 2016). We refer the reader for the work of Beirlant, Dudewicz, Gy\"oria and van der Meulen (1997) for a comprehensive review for nonparametric entropy estimators.

On the other hand, Bayesian estimation of entropy has not received much attention.  Exceptions include the work of  Mazzuchi, Soofi and Soyer (2000, 2008), who develop a Bayes estimate of $H(P)$ based on the Dirichlet process (Furguson, 1973) and provided a computational algorithm for their procedure.  The main goal of this paper is to derive an efficient and easy-to-implement Bayesian nonparametric estimator of (\ref{entropy}).  The anticipated estimator may be viewed as the Bayesian nonparametric counterpart  of the estimator  of Ebrahimi, Pflughoeft and Soofi (1994). A main motive of having this estimator, among others,  is  to use it in Bayesian methods such as model checking as discussed, for instance, in Al-Labadi and Evans (2018).  Therefore, it can be worthwhile  to have such an estimator in practice.

The remainder of this paper is organized as follows. In Section 2, the Dirichlet process prior is briefly reviewed. In Section 3, a Bayesian non-parametric estimator of the entropy is obtained and several of its properties are derived.  Section 4 develops a computational algorithm of the approach, where particular choices of $m$ and the hyperparameters of the Dirichlet process should be used.  Section 5 presents a number of examples where the behavior of the estimator is examined in some detail. A comparison between the new estimator,  the estimator of Vasicek's estimator and the estimator of Ebrahimi, Pflughoeft and Soofi  is also considered. Section 6 ends with a brief summary of the results.  Proofs are placed in the Appendix.

\subsection{Dirichlet process} \label{Dirichlet}

A relevant summary of the Dirichlet process is presented in this section. The Dirichlet process, formally introduced in Ferguson (1973), is considered the most well-known and widely used prior in Bayesian nonparameteric inference. Let  $\mathfrak{X}$ be a space and $\mathcal{A}$ be  a $\sigma-$algebra of subsets of $\mathfrak{X}$. Let $G$ be a fixed probability measure on $(\mathfrak{X},\mathcal{A})$, called the \emph{base measure},  %
 and $a$ be a positive number, called the \emph{concentration parameter}. Following Ferguson (1973), a
random probability measure $P=\left\{  P(A)\right\}  _{A\in\mathcal{A}}$ is
called a Dirichlet process on $(\mathfrak{X},\mathcal{A})$ with parameters $a$
and $G$, denoted by $DP(a,G)$, if for any finite measurable partition $\{A_{1},\ldots,A_{k}\}$ of
$\mathfrak{X}$ with $k \ge 2$, $\left(  P(A_{1}%
),\ldots\,P(A_{k})\right)\sim \text{Dirichlet}(aG(A_{1}),\ldots,$ $aG(A_{k}))$. It is assumed that if
$G(A_{j})=0$, then $P(A_{j})=0$ with a probability one.  For any
$A\in\mathcal{A},$ $P(A) \sim \text{Beta}(aG(A),a(1-G(A))$ and so ${E}(P(A))=G(A)\ \ $ and ${Var}(P(A))=G(A)(1-G(A))/(1+a).$
Thus, $G$ can be viewed as the center of the process. On the other hand, $a$ controls concentration, as the larger value of $a$, the more likely that $P$ will be close to $G$.

An important feature of the Dirichlet process is the conjugacy property. Specifically, if
$x=(x_{1},\ldots,x_{n})$ is a sample from $P\sim DP(a,G)$, then the posterior
distribution of $P$ is $P\,|\,x=P_x\sim DP(a+n,G_{x})$ where
\begin{equation}
G_{x}=a(a+n)^{-1}G+n(a+n)^{-1}F_{n}, \label{DP_posterior}%
\end{equation}
with $F_{n}=n^{-1}\sum_{i=1}^{n}\delta_{{x}_{i}}$ and $\delta_{x_{i}}$ the
Dirac measure at $x_{i}.$ Notice that, $G_{x}$
is a convex combination of the prior base distribution and the empirical
distribution. Clearly,   $G_{x}\to G$  as $a \to \infty$ while   $G_{x}\to F_n$ as $a \to 0$. On the other hand, by Glivenko-Cantelli theorem, when $a/n \to 0$ (i.e., $a$ is small comparable to $n$),  $G_{x}$ converges to true distribution function.  We refer the reader to Al-Labadi and Zarepour (2013a,b; 2014a) and Al-Labadi and Abdelrazeq (2017) for other asymptotic properties of the Dirichlet process.

Following Ferguson (1973), $P\sim{DP}(a,G)$ has the following series representation
\begin{equation}
P=\sum_{i=1}^{\infty}J_{i} \delta_{Y_{i}}, \label{series-dp}%
\end{equation}
where $\Gamma_{i}=E_{1}+\cdots+E_{i}$, $E_{i} \overset{i.i.d.}\sim \text{exponential}(1)$, $Y_{i} \overset{i.i.d.}\sim  G$
 independent of $\Gamma_{i}$, $L(x)=a\int_{x}^{\infty}t^{-1}e^{-t}dt,x>0,$ $L^{-1}(y)=\inf
\{x>0:L(x)\geq y\}$ and $J_{i}=L^{-1}(\Gamma_{i})/\sum_{i=1}^{\infty
}{L^{-1}(\Gamma_{i})}$.  It follows clearly from (\ref{series-dp}) that a
realization of the Dirichlet process is a discrete probability measure.  This
is correct even when $G$ is absolutely continuous. We refer to  $(Y_{i})_{i \ge 1}$ and  $(J_{i})_{i \ge 1}$  as the \emph{atoms} and the \emph{weights}, respectively.  Note that, one could resemble the discreteness of $P$ with the discreteness of $F_n$. Since data is always measured to finite accuracy, the true distribution being sampled from is discrete. This makes the  discreteness property of $P$ with no practical significant limitation.  

Because there is no closed form for the inverse of L\'evy measure $L(x)$, using Ferguson (1973) representation of the Dirichlet process is difficult in practice. As an alternative, Sethuraman (1994) uses  the stick-breaking approach to define the Dirichlet Process. Let $(\beta_i)_{i \ge 1}$ be a sequence of i.i.d. random
variables with a $\text{Beta}(1,\alpha)$ distribution. In (\ref{series-dp}), set
\begin{equation}
J_1=\beta_1,~ J_i=\beta_i\prod^{i-1}_{k=1}(1-\beta_k),~i\ge2.\label{eq6}
\end{equation}
and $(Y_{i})_{i \ge {1}}$ independent of $(\beta_i)_{i \ge 1}$. Unlike Ferguson's approach, the stick-breaking construction does not need normalization. By truncating the higher order terms in the sum to  simulate Dirichlet process, we can approximate the Sethuraman stick breaking representation by
\begin{equation*}
P_N=\sum_{i=1}^{N}J_{i,N} \delta_{{Y}_i}. \label{suth}
\end{equation*}
In here,  $(\beta_i)_{i \ge 1}$, $(J_{i,N})_{i \ge 1}$, and $(Y_i)_{i \ge 1}$ are as defined in (\ref{eq6}) with $\beta_N=1$. The assumption that $\beta_N=1$ is necessary to   make the weights add up to 1 almost surely (Ishwaran and James, 2001).

The Dirichlet process can also be obtained from the following finite mixture
models developed by Ishwaran and Zarepour (2002). Let $P_{N}$ has the from given in (\ref{series-dp}) with $(J_{1,N}, \ldots, J_{N,N})\sim\,$Dirichlet$(a/N,\ldots,a/N)$.
Then $E_{P_{N}}(g)\rightarrow E_{P}(g)$ in distribution as
$N\rightarrow\infty$, for any measurable function $g:\mathbb{R}\rightarrow
\mathbb{R}$ with $\int_{\mathbb{R}}|g(x)|H(dx)<\infty$ and $P\sim DP(a,G)$. In particular, $(P_{N})_{N\geq1}$ converges in distribution to $P$, where $P_{N}$ and $P$ are random values in the space $M_{1}(\mathbb{R})$ of probability measures on $\mathbb{R}$ endowed with the topology of weak convergence. To generate $(J_{i,N})_{1\leq i\leq N}$ put $J_{i,N}=G_{i,N}%
/\sum_{i=1}^{N}{G_{i,N},}$ where $(G_{i,N})_{1\leq i\leq N}$ is a sequence of i.i.d. gamma$(a/N,1)$ random variables independent of $(Y_{i})_{1\leq i\leq N}$.

For other simulation methods for the Dirichlet process, see
Bondesson (1982),  Wolpert and Ickstadt (1998) and  Zarepour and Al-Labadi (2012), Al-Labadi and Zarepour (2014b).

\section{Bayesian Estimation of the Entropy}
\label{sec-entropy}
Let $P_N=\sum_{i=1}^{N}J_{i,N} \delta_{{Y}_i}$.  Similar to (\ref{slope}),
the slope of the straight line that joins the two points $\left(P_N(Y_{i-m})= \sum_{k=1}^{i-m} J_{k,N},
~Y_{(i-m)}\right)$ and $\Big(P_N(Y_{(i+m)})=\sum_{k=1}^{i+m} J_{k,N},$ $~Y_{(i+m)}\Big)$ is
\begin{align*}
\frac{Y_{(i+m)}-Y_{(i-m)}}{P_N(Y_{(i+m)})-P_N(Y_{(i-m)})}=\frac{Y_{(i+m)}-Y_{(i-m)}}{c_{i,a}},
\end{align*}
where
\begin{align}
c_{i,a} = \left\{
        \begin{array}{ll}
          \sum_{k=2}^{i+m}J_{k,N} & \quad  1\leq i \leq m \\
           \sum_{k=i-m+1}^{i+m}J_{k,N} & \quad  m+1 \leq i \leq N-m \label{coef}\\
            \sum_{k=i-m+1}^{N}J_{k,N} & \quad  N-m+1 \leq i \leq N
        \end{array}
    \right..
\end{align}
Note that,  from the properties of the Dirichlet distribution, we have $J_{i,N}$ $\sim Beta\left(a/N,a(1-1/N)\right)$. Thus, $E\left(  J_{i,N}\right)
=N^{-1}$. Hence,
\begin{align*}
E[c_{i,a}] &= \left\{
        \begin{array}{ll}
            \frac{m+i-1}{N} & \quad  1\leq i \leq m \\
            \frac{2m}{N} & \quad  m+1 \leq i \leq N-m\\
           \frac{N+m-i}{N} & \quad  N-m+1 \leq i \leq N\\
        \end{array}
    \right.\\
    &=\frac{m}{N}c_{i},
\end{align*}
where $c_i$ is defined in  (\ref{EPS1}). The next proposition underlines a direct connection  between $c_{i,a}$ and $c_{i}$. Its proof is given in the Appendix.
\begin{proposition} \label{lemma-expectation0} Let $(J_{i,N})_{1\leq i\leq N}\sim\,$Dirichlet$(a/N,\ldots,a/N)$. As $N \to \infty$,
\begin{enumerate}
  \item $J_{i,N}-1/N \overset{p} \to 0$
  \item $c_{i,a}-\frac{m}{N}c_{i} \overset{p} \to 0$, where $c_{i,a}$ and $c_{i}$ are defined  in (\ref{coef}) and (\ref{EPS}), respectively.
\end{enumerate}
\end{proposition}

Proposition  \ref{lemma-expectation0} motivates the possibility of constructing a Bayesian non-parametric version of (\ref{EPS1}) based on the Dirichlet process. The precise form of the anticipated estimator to (\ref{entropy}) is presented in the next lemma. The proof is placed in the Appendix.

\begin{lemma} \label{lemma-expectation} Let $P_N=\sum_{i=1}^{N}J_{i,N} \delta_{{Y}_i}$ as defined in Section \ref{Dirichlet}, where  $Y_{1},Y_{2}, $ $\ldots, Y_{N}\overset{i.i.d.}\sim G$. Let $m$ be a positive integer smaller than $N/2$, $Y_{(i)}=Y_{(1)}$ if $i<1$, $Y_{(i)}=Y_{(N)}$ if $i>N$ and $Y_{(1)} \le Y_{(2)} \le \cdots \le Y_{(N)}$ are the order statistics of $Y_{1}, Y_{2}, \ldots, Y_{N}$. Let
\begin{align}
H_{m,N,a}=\frac{1}{N}\sum_{i=1}^N  \log \left(\frac{Y_{(i+m)}-Y_{(i-m)}}{c_{i,a}}\right), \label{entropy_Bayesian}
\end{align}
where $c_{i,a}$ is defined in (\ref{coef}). As $N \to \infty$, $m \to \infty$, $m/N \to 0$ and $a \to \infty$, we have
\begin{align*}
E\left[H_{m,N,a}\right]-E\left[H^{EPS}_{m,n}\right] \to 0,
\end{align*}
where $H^{EPS}_{m,n}$ is defined in (\ref{EPS}).
\end{lemma}

The next lemma shows that the estimator defined in (\ref{entropy_Bayesian})  is consistent. The defined a Bayesian nonparametric prior for the entropy. A formal proof is given in the Appendix.

\begin{lemma} \label{lemma-entropy_prior} Let  $H_{m,N,a}, N, m, a$ and $G$ be as defined in Lemma \ref{lemma-expectation}. Then
as $N \to \infty$, $m \to \infty$, $m/N \to 0$ and $a \to \infty$, we have
\begin{align*}
H_{m,N,a} \overset{p} \to H(G)=-\int_{-\infty}^{\infty} g(x)\log g(x) dx,
\end{align*}
where $G^{\prime}(x)=g(x)$.
\end{lemma}

\section{Computations and the Choices of $m$, $a$ and $G$}
Let  $x=(x_{1},\ldots,x_{n})$ be a sample from a continuous distribution $P$. The aim
is to approximate $H(P)$ as defined in (\ref{entropy}). We will use the prior $P\sim
DP(a,G)$ for some choice of $a$ and $G$  so $P\,|\,x\sim DP\left(
a+n,G_{x}\right)$. See Section 2.

To fully implement the approximation $H_{m,N,a}$ as in Lemma \ref{lemma-entropy_prior}, it is necessary to discuss the choices for $m$, $a$ and $G$. We start by the choice of $m$, where its optimal value   is still an open problem in entropy estimation. However, as discussed in Vasicek (l976), with increasing $N$, the best value of $m$ increases while the ratio $m/N$ tends to zero. For example, for $N=10, 20, 50$ Vasicek  (1976) recommended using $m=2, 3, 4$, respectively.  On the other hand, Grzegorzewski and Wieczorkowski (1999) proposed the following formula for optimal values of $m$:
 \begin{align}
\label{optimal-m} m =\floor {\sqrt{N} +0.5},
    \end{align}
where $\floor {y}$ is the largest integer less than or equal to $y$. Thus, by (\ref{optimal-m}), for $N=10, 20, 50$, the best choices of $m$ are 2, 3, 7,  respectively. In this paper, we will use the rule  (\ref{optimal-m}).
Note that, the value of $m$ in (\ref{optimal-m}) is the value that will be used for the prior. For the posterior, one should replace $N$ by the number of distinct atoms in $P_N|\,x$,  an approximation of $P|x$. Observe that, it follows from (\ref{DP_posterior}) that if $a/n$ is close to zero, then the number of distinct  atoms in $P_N|\,x$ will typically be  $n$.

As for hyperparameters $a$ and $G$, their choices depend on the application of interest. For instance, for model checking,  to detect small deviations, $a$ and $G$ should be selected so that there is a good concentration about the prior (Al-Labadi and Evans, 2018). Further, they  recommended that $a$ should be chosen so that its  value does not exceed  $0.5n$ as otherwise the prior may become too influential. In light of this under the context of entropy estimation, any choice of $a$ such that $a/n$ is close to zero should be compatible with any choice of $G$. This follows from  (\ref{DP_posterior}) as the sample will dominate the prior guess $G$.  For example, setting $a=0.05$ and $n=10$ in (\ref{DP_posterior}) gives
\begin{equation*}
G_{x}=0.005G+0.995F_{n},
\end{equation*}
which means the chance to draw a sample from the collected data is 99.5\% over a new sample from $G$. For simplicity, we suggest to set $G=N(0,1)$ and $a=0.05$, although other choices are certainly possible.  An example studying the  sensitivity of the approach to the choice of $G$ is covered in Section 4.

The following result shows that, as the sample size increases (i.e., the concentration parameter $a$ is small comparable to the sample size $n$),  then the posterior of $H_{m,N,a}$  (i.e., the proposed estimator) converges in probability to (\ref{entropy}). The proof follows from (\ref{DP_posterior}), Glivenko-Cantelli theorem and Lemma 3.

\begin{lemma} \label{lemma-entropy_post} Let $x=(x_{1},\ldots,x_{n})$ be a sample from $P\sim DP(a,G)$. Let $H_{m,N,a}$ be as defined in Lemma  \ref{lemma-expectation}. Then
as $N \to \infty$, $m \to \infty$, $n \to \infty$, $m/N \to 0$ and $a/n \to 0$
\begin{align*}
H_{m,N,a}|x \overset{p} \to H(P)=-\int_{-\infty}^{\infty} p(x)\log p(x) dx.
\end{align*}
\end{lemma}

Now, based on Lemma \ref{lemma-entropy_post}, the following gives a computational algorithm for estimating (\ref{entropy}).

\vspace{0.3 cm}

\noindent\textbf{Algorithm A}\textit{(}\emph{Nonparametric Estimation of Entropy}\textit{):}
\begin{enumerate}
\item [(i)] Let $P \sim DP(a,G)$ and $P_N$ be an approximation of $P$.  Set $a=0.05$ and $G=N(0,1)$.
\item [(ii)]  Generate  a sample from $P_N|x$, where $P_N|x$ is an approximation of   $P|x \sim DP(a+n,G_{x})$. See Section 2.
\item [(iii)] Compute $H_{m,N,a}|x$ as in Lemma \ref{lemma-entropy_post}.
\item [(iv)]  Repeat steps (i) and (iii) to obtain a sample of $r$ values from $H_{m,N,a}|x$. For large $r$, the empirical distribution of these values is an approximation to the distribution of $H_{m,N,a}|x$.
\item  [(v)] The average of the $r$ values  generated in step (iv) will be the estimator of the entropy.
\end{enumerate}

Note that, for estimation purposes,  the prior has no significant role. This is not necessarily will be the case for other applications such as model checking.

\section{Examples}

In this section, we study the behaviour of the proposed  estimator  in terms of efficiency and robustness. The proposed estimator is   compared with its non-Bayesian counterpart estimators of Vasicek  (1976) and Ebrahimi, Pflughoeft and Soofi (1994). Additionally, for a comprehensive comparison, we included the (weighted) Kozachenko–Leonenko entropy estimator (Kozachenko and Leonenko, 1987; Berrett, Samworth and Yuan, 2018), which is based on the $k$-nearest neighbour distances of the sample. We set the value of $k$ to equal to $m$ in (\ref{optimal-m}). This value of $k$ (square root of the sample size) is recommended, for instance, by Mitra, Murthy and Pal (2002) and  Bhattacharya, Ghosh,  and Chowdhur (2012). For each sample size ($n=10, 20, 50$), 1000 samples were generated. We have considered four distributions: exponential with mean 1 (exact entropy is 1), Uniform  on $(0,1)$ (exact entropy is 0), $N(0,1)$ (exact entropy is $0.5\log(2\pi e)$) and Weibull distribution with shape parameter equal to 2 and scale parameter equal to 0.5 (exact value of $-0.0977$).  The estimators and  their mean squared  errors are computed. Here each sample of the 1000 samples gives an estimate. The reported value of the estimator (Est) is the average of  the 1000 estimates. On the other hand, the mean squared error (MSE) is computed as follows: $(\mbox{estimated value for each sample}- \mbox{true value})^2/1000$. The computing program codes were implemented in the programming language \textbf{\textsf{R}}. For the KL entropy estimator, we used the  package \textbf{IndepTest} (Berrett, Grose and Samworth, 2018). In all  cases, the prior was taken to be $DP\left(a,N(0,1)\right)$. In Algorithm A,  we set $r=1000$ and $N=200$. The sensitivity to the choice of $\ a$  is  investigated  and we record only a few values in the tables.

\begin{table}[htbp]
\caption{Uniform $(0,1)$.} \label{tab3}
  \centering
      \begin{tabular}[c]{lllccccccc}
    \hline
 & & &\multicolumn{1}{c}{$H_{m,N,a}|x$}  &\multicolumn{1}{c}{$H^{V}_{m,n,a}$}  &\multicolumn{1}{c}{$H^{EPS}_{m,n,a}$} & \multicolumn{1}{c}{KL Entropy} \\ \cline{4-4} \cline{5-5} \cline{6-6}  \cline{7-7} $n$&$m$ &$a$&Est(MSE)  & Est(MSE) & Est(MSE)&st(MSE)  \\
    \hline
10&3&0.05  &  $-0.017(0.017)$ & $-0.410(0.193)$ &  $-0.154(0.048)$ & $0.073(0.080)$\\
& & 1  &  0.517(0.275)& &  &\\
& & 5  & 0.767(0.593) &   &  & \\
\hline
20	&4& 0.05  & $-0.050(0.010)$ & $-0.260(0.077)$ & $-0.102(0.017)$ & 0.038(0.035)\\
	& & 1  & 0.418(0.179) &  &  & \\
	& & 5  &  0.657(0.435)&  &  &  \\
\hline
50	&7& 0.05  &$-0.051(0.004)$ & $-0.150(0.024)$  & $-0.051(0.004)$ & $0.014(0.013)$ \\	
	& & 1  & $0.230(0.055)$  &&& \\
	& & 5  & $0.496(0.247)$  &&& \\
\hline
    \end{tabular}%
\end{table}%

\begin{table}[htbp]
\caption{Exponential with mean 1.} \label{tab3}
  \centering
      \begin{tabular}[c]{lllccccc}
 & & &\multicolumn{1}{c}{$H_{m,N,a}|x$}  &\multicolumn{1}{c}{$H^{V}_{m,n,a}$}  &\multicolumn{1}{c}{$H^{EPS}_{m,n,a}$} & \multicolumn{1}{c}{KL Entropy} \\ \cline{4-4} \cline{5-5} \cline{6-6}  \cline{7-7} $n$&$m$ &$a$&Est(MSE)  & Est(MSE) & Est(MSE)&Est(MSE)  \\
    \hline
10	&3& 0.05  &  0.891(0.114) & 0.575( 0.298) &  0.831(0.146) & 0.947(0.149)  \\
&& 1  & 1.057(0.073) &&& \\
&& 5  & 1.175(0.068)  &&&  \\
\hline
20	&4& 0.05  &0.937(0.052) & 0.752(0.112) &  0.910(0.059) & 0.967(0.069)\\
	&& 1  & 1.148(0.066)& &&\\
	&& 5  & 1.208(0.075) &&& \\
\hline
50	&7& 0.05 &  0.956(0.022) & 0.864(0.039)  & 0.964(0.022) & 0.967(0.026)  \\
&& 1  & 1.137(0.039)& && \\
	&& 5  & 1.222(0.068)& && \\
\hline
    \end{tabular}%
\end{table}%

\begin{table}[htbp]
\caption{$N(0,1)$.} \label{tab3}
  \centering
      \begin{tabular}[c]{lllccccc}
 & & &\multicolumn{1}{c}{$H_{m,N,a}|x$}  &\multicolumn{1}{c}{$H^{V}_{m,n,a}$}  &\multicolumn{1}{c}{$H^{EPS}_{m,n,a}$} & \multicolumn{1}{c}{KL Entropy} \\ \cline{4-4} \cline{5-5} \cline{6-6}  \cline{7-7} $n$&$m$ &$a$&Est(MSE)  & Est(MSE) & Est(MSE)&Est(MSE)  \\
    \hline
10	&3& 0.05  &1.112(0.159)& 0.869(0.374) &    1.1253(0.158)&   1.293(0.124)\\
	&& 1  &1.183(0.087)  & & & \\
	&& 5  &1.207(0.061)  & & & \\
\hline
20	&4& 0.05 & 1.223(0.069) & 1.092(0.138)  &  1.251(0.060) &  1.344(0.058)\\
	&& 1  &1.251(0.049)   &&& \\
	&& 5  &1.282(0.031)   &&&\\
\hline
50	&7& 0.05  & 1.331(0.020) & 1.258(0.038) &   1.358(0.016)  & 1.388(0.022) \\
	&& 1  &1.332(0.018)  &&&   \\
	&& 5  &1.340(0.014)  &&&   \\
\hline
    \end{tabular}%
\end{table}%

\begin{table}[htbp]
\caption{Weibull with shape parameter 2 and scale parameter 0.5.} \label{tab3}
  \centering
      \begin{tabular}[c]{lllccccc}
 & & &\multicolumn{1}{c}{$H_{m,N,a}|x$}  &\multicolumn{1}{c}{$H^{V}_{m,n,a}$}  &\multicolumn{1}{c}{$H^{EPS}_{m,n,a}$} & \multicolumn{1}{c}{KL Entropy} \\ \cline{4-4} \cline{5-5} \cline{6-6}  \cline{7-7} $n$&$m$ &$a$&Est(MSE)  & Est(MSE) & Est(MSE)&Est(MSE)  \\
    \hline
10	&3& 0.05  & $-0.200(0.056)$ & $-0.641(0.363)$ & $-0.385(0.150)$ &$-0.209(0.123)$ \\
	&& 1  &$0.379(0.242)$  &&&  \\
	&& 5  &$0.652(0.572)$  &&&  \\
\hline
20	&4& 0.05  & $-0.195(0.033)$& $-0.417(0.129)$ & $-0.258(0.053)$ &$-0.148(0.055)$\\
	&& 1  &$0.278(0.150)$  &&&  \\
	&& 5  &$0.509(0.376)$  &&&  \\
\hline
50	&7& 0.05  &  $-0.170(0.016)$ & $-0.269(0.040)$ & $-0.169(0.016)$& $-0.131(0.023)$\\
	&& 1  &$0.094(0.043)$ &&&  \\
	&& 5  &$0.330(0.187)$  &&& \\
\hline
    \end{tabular}%
\end{table}%

It follows clearly from Table 1 - Table 4 that, when $a=0.05$, the new approximation of entropy has the lowest mean squared error for most cases.  As illustrated in Section 3, the choice of $a$ is extremely important and for the case of estimation it should be chosen so that $a/n$ is close to zero. The choice $a=0.05$ is found to be satisfactory in all the cases considered in the paper.

It is also interesting to consider the effect of using different base measures $G$ on the
methodology. We fix $a$ at 0.05 and 5. We used several values of $G$. To this end, the next data set is generated from the exponential distribution with mean 20.
\begin{quote}
$1.884,5.289, 20.890,  20.093,  21.007,  15.261,   7.716,
18.979,  27.537, 10.291,    \newline%
 31.048, 1.215,  13.564,  14.966, 24.896,  10.849$
\end{quote}

The results of the estimated entropy for the previous data set are  reported in Table 5. Clearly, using different $G$ with $a=0.05$ has no impact on the estimated value. However, when $a=5$, the estimated value depends on the choice of $G$.

\begin{table}[htbp] \centering
\begin{tabular}
[c]{ccc}\hline
$G$  & Estimate: $a=0.05$ & Estimate: $a=5$ \\
\hline
$N(0,1)$&3.402 &3.352\\
$N(3,9)$&3.407&3.393\\
$t_{1}$ &3.407&3.735\\
${\cal E}(1))$&3.402&3.143\\
$U[0,1]$&3.398&3.118\\
\hline

\end{tabular}
\caption{Study of the effect of the proposed estimator using different base measures $G$ of $P \sim DP(a=0.05,G)$. Here $N(\mu,\sigma^2)$ is the normal distribution with mean $\mu$ and standard deviation $\sigma$, $t_1$ is the $t$ distribution with $1$ degrees of freedom, ${\cal E}(1)$ is the exponential distribution with mean $1$ and  $U[0,1]$ is the uniform distribution over $[0,1]$.}
\end{table}%

\section{Conclusion}
In this paper, an efficient yet simple Bayesian nonparametric estimator of entropy is proposed.  The proposed estimator is considered an analogous Bayesian estimator to the estimator of Ebrahimi, Pflughoeft and Soofi (1994). Through several examples, it has been shown that the approach performs extremely well where a smaller mean squared error is obatined. A foremost motive of having this estimator to use it in applications such as model checking as discussed, for instance, in Al-Labadi and Evans (2018). We have left this critical avenue of research to future work.

%
%

\renewcommand{\theequation}{\thesection.\arabic{equation}}
\setcounter{equation}{0}

\appendix

\section{Proofs}

\subsection{Proof of Proposition 1}

{1.} Note that, for each $1 \le i \le N$, from the properties of the Dirichlet distribution, $J_{i,N}\sim$B$eta\left(a/N,a(1-1/N)\right)$. It follows that $E\left[J_{i,N}\right]=1/N$ and
$$V\left[J_{i,N}\right]=\frac{1/N(1-1/N)}{a+1},$$
where $V$ stands for the variance. Since, as $N \to \infty$, $V\left(J_{i,N}\right)\to 0$,  we conclude the result.

\smallskip

{2.} By ($\ref{expectation}$), $E[c_{i,a}]=\frac{m}{N}c_{i}$. Thus,  to prove the proposition, by Chebyshev's inequality, it is sufficient to show that $V(c_{i,a}) \to 0$.  We consider three cases.

\smallskip\underline{Case I (for $1 \leq i \leq m$):} From the aggregation property of the Dirichlet distribution,
\begin{align}
 \sum_{k=2}^{i+m}J_{k,N}\sim B{\rm{eta}}\left(\sum_{k=2}^{i+m} \frac{a}{N},a-\sum_{k=2}^{i+m} \frac{a}{N}\right). \label{sum1}
\end{align}
Hence,
\begin{align*}
V\left(\sum_{k=2}^{i+m}J_{k,N}\right)&=\frac{\sum_{k=2}^{i+m} \frac{a}{N}\left(a-\sum_{k=2}^{i+m} \frac{a}{N}\right)}{a^2(1+a)}\\
&= \frac{(i+m-1)(N-i-m+1)}{N^2(a+1)} \to 0,
\end{align*}
as $N \to \infty$.

\smallskip

\underline{Case II (for $m+1 \leq i \leq N-m$):} similar to Case I,
\begin{align}
\sum_{k=i-m+1}^{i+m}J_{k,N}\sim B{\rm{eta}}\left(\sum_{k=i-m+1}^{i+m} \frac{a}{N},a-\sum_{k=i-m+1}^{i+m} \frac{a}{N}\right). \label{sum2}
\end{align}
Hence,
\begin{align*}
V\left(\sum_{k=i-m+1}^{i+m}J_{k,N}\right)&=\frac{\sum_{k=i-m+1}^{i+m} \frac{a}{N}\left(a-\sum_{k=i-m+1}^{i+m} \frac{a}{N}\right)}{a^2(a+1)} \\
&=\frac{2m(N-2m)}{N^2(a+1)}  \to 0,
\end{align*}
as $N \to \infty$.

\smallskip

\underline{Case III (for $N-m+1 \leq i \leq N$):} As in the previous cases,
\begin{align}
\label{sum3} \sum_{k=i-m+1}^{N}J_{k,N}\sim B{\rm{eta}}\left(\sum_{k=i-m+1}^{N} \frac{a}{N},a-\sum_{k=i-m+1}^{N} \frac{a}{N}\right).
\end{align}
Therefore,
\begin{align*}
V\left(\sum_{k=i-m+1}^{N}J_{k,N}\right)&=\frac{\sum_{k=i-m+1}^{N} \frac{a}{N}\left(a-\sum_{k=i-m+1}^{N} \frac{a}{N}\right)}{a^2(1+a)}\\
&= \frac{(N-i+m)(i-m)}{N^2(a+1)}\to 0,
\end{align*}
as $N \to \infty$. This complete the proof of the proposition. \endproof

\bigskip
\subsection{Proof of Lemma 2}
Recall that,
\begin{align*}
H_{m,N,a}&=\frac{1}{N}\sum_{i=1}^N  \log \left(\frac{Y_{(i+m)}-Y_{(i-m)}}{c_{i,a}}\right)
\end{align*}
and
\begin{align*}
H^{EPS}_{m,N}=&\frac{1}{N}\sum_{i=1}^N  \log \left(\frac{Y_{(i+m)}-Y_{(i-m)}}{mc_i/N}\right),
\end{align*}
where $c_{i,a}$ and $c_{i}$ are defined, respectively, on (\ref{coef}) and (\ref{EPS1}). Thus,
        \begin{align}
      \nonumber E\left[H_{m,N,a}\right]-E\left[H^{EPS}_{m,N}\right] & =\frac{1}{N}\sum_{i=1}^N E\left[\log \left(\frac{mc_i/N}{c_{i,a}}\right)\right]\\
       \label{con3} & =\frac{1}{N}\sum_{i=1}^N \log \left(c_i m/N\right)-\frac{1}{N}\sum_{i=1}^N E\left[\log c_{i,a}\right].
        \end{align}
We want to show that, as $N \to \infty$, $m \to \infty$, $m/N \to 0$ and $a \to \infty$,  (\ref{con3}) $\to 0$. We consider three cases.
\smallskip

\underline{Case I (for $1 \leq i \leq m$):} notice that, with  $\alpha_k=aN^{-1}$ and $\alpha_0=\sum_{k=1}^N \alpha_k=a$, $E\left[J_{i,N} \log c_{i,a}\right]=E\left[J_{i,N}\log \left(\sum_{k=2}^{i+m}J_{k,N}\right)\right]=$
\begin{align}
\nonumber &\int \cdots \int z_{i}\log \left(\sum_{k=2}^{i+m} z_{k}\right) \frac{\Gamma(\alpha_0)}{\prod_{k=1}^N \Gamma(\alpha_k)}z_1^{\alpha_1-1}\cdots z^{\alpha_N-1}  dz_1 \cdots dz_i\cdots dz_N, \\
\nonumber &=\int \cdots \int \log \left(\sum_{k=2}^{i+m} z_{k}\right) \frac{\Gamma(\alpha_0)}{\prod_{k=1}^N \Gamma(\alpha_k)}z_1^{\alpha_1-1}\cdots z_i^{\left(\alpha_i+1\right)-1} \cdots z^{\alpha_N-1} \\
\nonumber &dz_1 \cdots dz_i\cdots dz_N \\
\label{convergence2-caseI} &= \frac{\alpha_i}{\alpha_0}E\left[\log \left(\sum_{k=2}^{i+m}Z_{k,N}\right)\right],
\end{align}
where $\sum_{k=2}^{i+m}Z_{k,N}\sim Beta\left(\sum_{k=2}^{i+m}\alpha_{k}+1, (a+1)-\left(\sum_{k=2}^{i+m}\alpha_{k}+1\right)\right)$. For $\alpha_k=aN^{-1}$, $\sum_{k=2}^{i+m}Z_{k,N}\sim Beta\left(a(m+i-1){N}^{-1}+1, a-a(m+i-1){N}^{-1}\right)$. From the properties of the beta distribution, we have
\begin{align}
 \nonumber (\ref{convergence2-caseI})&=\frac{\alpha_i}{\alpha_0}\left(\psi\left(\sum_{k=1}^{m}\alpha_{k}+1\right)-\psi\left(\alpha_0+1\right)\right)\\
 &=\frac{1}{N} \bigg(\psi\left(\frac{a(m+i-1)}{N}+1\right)-\psi(a+1)\bigg), \label{convergence-aa}
\end{align}
where $\psi(x) = \Gamma'(x)/\ \Gamma(x)$ is the digamma function. Therefore, by (\ref{convergence-aa}) and for $1 \leq i \leq m$, we obtain
\begin{align*}
(\ref{con3})&=\frac{1}{N}\sum_{i=1}^m\log \left(\frac{m+i-1}{N}\right)\\
&-\frac{1}{N}\sum_{i=1}^m\left(\psi\left(\frac{a(m+i-1)}{N}+1\right)-\psi(a+1)\right).
\end{align*}
Using that facts that $\psi (x+1)=\log (x)+O\left(x^{-1}\right)$ and $\sum_{i=0}^{L-1}\frac{1}{x+i}=\psi(x+L)-\psi(x)=\log\left(\frac{x+L}{x}\right)+O\left(x^{-1}\right)$ (Abramowitz and Stegun, 1972), we have
\begin{align*}
(\ref{con3})&=-\frac{1}{N}\sum_{i=1}^m O\left(\frac{N}{a(m+i-1)}\right)+\frac{1}{N}\sum_{i=1}^m O\left(\frac{1}{a}\right)\\
&=O\left(\sum_{i=1}^m\frac{1}{a(m+i-1)}\right)+ O\left(\frac{m}{Na}\right)\\
&=\frac{1}{a}O\left(\psi(2m)-\psi(m-1)\right)+ O\left(\frac{m}{Na}\right)\\
&=\frac{1}{a}O\left(\log\left(\frac{2m}{m-1}\right)+\frac{1}{2m}\right)+  O\left(\frac{m}{Na}\right) \to 0
\end{align*}

as $N \to \infty$, $m \to \infty$, $m/N \to 0$ and $a \to \infty$.

\smallskip

\underline{Case II (for $m+1 \leq i \leq N-m$):} similar to Case I,
\begin{align}
E\left[J_{i,N} \log c_{i,a}\right]&= \frac{\alpha_i}{\alpha_0}E\left[\log \left(\sum_{k=i-m+1}^{i+m}Z_{k,N}\right)\right],\label{convergence2-caseII}
\end{align}
where $\sum_{k=i-m+1}^{i+m}Z_{k,N}\sim Beta\left(2am{N}^{-1}+1, a-2am{N}^{-1}\right)$. Thus, from the properties of the beta distribution, we have
\begin{align}
(\ref{convergence2-caseII})=\frac{1}{N} \left(\psi\left(\frac{2am}{N}+1\right)-\psi(a+1)\right). \label{convergence-b}
\end{align}
Therefore, by (\ref{convergence-b}), we have
\begin{align*}
(\ref{con3})&=\sum_{i=m+1}^{N-m}\log \left(\frac{2m}{N}\right)-\frac{1}{N}\sum_{i=m+1}^{N-m}\left(\psi\left(\frac{2am}{N}+1\right)+\psi(a+1)\right)\\
&=O\left(\frac {N}{2am}\right)+O\left(\frac{1}{a}\right) \to 0
\end{align*}
as $N \to \infty$, $m \to \infty$, $m/N \to 0$ and $a \to \infty$.

\smallskip

\underline{Case III (for $N-m+1 \leq i \leq N$):} similar to the previous cases,
\begin{align}
E\left[J_{i,N} \log c_{i,a}\right]= \frac{\alpha_i}{\alpha_0}E\left[\log \left(\sum_{k=i-m+1}^{N}Z_{k,N}\right)\right], \label{convergence2-caseIII}
\end{align}
where $\sum_{k=i-m+1}^{N}Z_{k,N}\sim Beta\left(a(N+m-i){N}^{-1}+1, a-a(N+m-i){N}^{-1}\right)$. Thus, from the properties of the beta distribution, we have
\begin{align}
(\ref{convergence2-caseIII})=\frac{1}{N} \left(\psi\left(\frac{a(N+m-i)}{N}+1\right)-\psi(a+1)\right). \label{convergence-c}
\end{align}
Therefore, by (\ref{convergence-c}), we have
\begin{align*}
(\ref{con3})&=\frac{1}{N}\sum_{i=N+m+1}^{N}\log \left(\frac{N+m-i}{N}\right)\\
&-\frac{1}{N}\sum_{i=N-m+1}^{N}\psi\left(\frac{a(N+m-i)}{N}+1\right)+\psi(a+1)\\
&= O\left(\sum_{i=N-m+1}^{N}\frac {1}{a(N+m-i)}\right)+O\left(\frac{m}{aN}\right)\\
&=O\left(\frac{\psi(1-2m)}{a}-\frac{\psi(1-m)}{a}\right)+O\left(\frac{m}{aN}\right)\\
&\to 0
\end{align*}
as $N \to \infty$, $m \to \infty$, $m/N \to 0$ and $a \to \infty$. Thus, in all cases, as $N \to \infty$, $m \to \infty$, $m/N \to 0$ and $a \to \infty$,  (\ref{con3}) $\to 0$.  This complete the proof of Lemma 2. \endproof

\subsection{Proof of Lemma 3}  Note that,
\begin{align*}
H_{m,N,a}=\left(H_{m,N,a}-H^{EPS}_{m,N}\right)+H^{EPS}_{m,N},
\end{align*}
where  $H^{EPS}_{m,N}$ is the approximation given in (\ref{EPS}). It follows that,
\begin{align*}
 H_{m,N,a}-H^{EPS}_{m,N}&=\frac{1}{N}\sum_{i=1}^N \log \left(\frac{mc_i/N}{c_{i,a}}\right).
\end{align*}
Since  $\left(J_{i,N}\right)_{1\le i \le N}$ is a sequence of  pairwise negative associated identically
distributed random variables with finite expectations, by Theorem 4.2.8 of Atkinson (2017), the weak law of large numbers holds for the sequence  $\left(\frac{mc_i/N}{c_{i,a}}\right)_{1\le i \le N}$. Thus we have
\begin{align*}
 H_{m,N,a}-H^{EPS}_{m,N} -E\left(\log \left(\frac{mc_i/N}{c_{i,a}}\right)\right) \to 0.
\end{align*}
We show that $E\left(\log \left(\frac{mc_i/N}{c_{i,a}}\right)\right) \to 0$. We consider three cases.

\smallskip\underline{Case I (for $1 \leq i \leq m$):} From (\ref{sum1}) and  the well-known property of the beta distribution, we have
\begin{align*}
E\left(\log \left(\frac{mc_i/N}{c_{i,a}}\right)\right)&= \log \left(\frac{mc_i}{N}\right)-E\left(\log \left(c_{i,a}\right)\right)\\
&=  \log \left(\frac{i+m-1}{N}\right) -\psi \left( \frac{a(i+m-1)}{N}\right)+\psi(a)\\
&=  \log \left(\frac{i+m-1}{N}\right) -\log\left(\frac{a(i+m-1)}{N}\right)-O\left(\frac{N}{a(i+m-1)}\right)\\
& +\log(a)+O(\frac{1}{a})\\
&= -O\left(\frac{N}{a(i+m-1)}\right)+O(\frac{1}{a}) \to 0.
\end{align*}
as  $N \to \infty$, $m \to \infty$, $m/N \to 0$ and $a \to \infty$.

\smallskip

\underline{Case II (for $m+1 \leq i \leq N-m$):} From (\ref{sum2}) and  the well-known property of the beta distribution, we have
\begin{align*}
E\left(\log \left(\frac{mc_i/N}{c_{i,a}}\right)\right)&= \log \left(\frac{mc_i}{N}\right)-E\left(\log \left(c_{i,a}\right)\right)\\
&=  \log \left(\frac{2m}{N}\right) -\psi \left( \frac{2am)}{N}\right)+\psi(a)\\
&=  \log \left(\frac{2m}{N}\right) -\log\left(\frac{2am}{N}\right)-O\left(\frac{N}{2am}\right)\\
& +\log(a)+O(\frac{1}{a})\\
&= -O\left(\frac{N}{2am}\right)+O(\frac{1}{a}) \to 0.
\end{align*}
as  $N \to \infty$, $m \to \infty$, $m/N \to 0$ and $a \to \infty$.

\smallskip

\underline{Case III (for $N-m+1 \leq i \leq N$):} As in the previous cases, from (\ref{sum3}), we have
\begin{align*}
E\left(\log \left(\frac{mc_i/N}{c_{i,a}}\right)\right)&= \log \left(\frac{mc_i}{N}\right)-E\left(\log \left(c_{i,a}\right)\right)\\
&=  \log \left(\frac{n+m-i}{N}\right) -\psi \left( \frac{a(n+m-i)}{N}\right)+\psi(a)\\
&=  \log \left(\frac{n+m-i}{N}\right) -\log\left(\frac{a(n+m-i)}{N}\right)-O\left(\frac{N}{a(n+m-i)}\right)\\
& +\log(a)+O(\frac{1}{a})\\
&= -O\left(\frac{N}{a(n+m-i)}\right)+O(\frac{1}{a}) \to 0.
\end{align*}
as  $N \to \infty$, $m \to \infty$, $m/N \to 0$ and $a \to \infty$. Thus, in all cases,  $H_{m,N,a}-H^{EPS}_{m,N} \overset{p} \to 0$. Also, by Ebrahimi, Pflughoeft and Soofi (1994), as $N \to \infty$, $m \to \infty$ and $m/N \to 0$ we have $H^{EPS}_{m,N}\overset{p}\to H(G)$. Now, applying  Slutsky's theorem (Furguson, 1996) completes the proof. \endproof


%


\begin{thebibliography}{99}                                                                                               %

\bibitem {bib1} Abramowitz, M., and Stegun, I. A. (1972). \emph{Handbook of Mathematical Functions with
Formulas, Graphs, and Mathematical Tables}. New York, Dover.

\bibitem {bib2} Alizadeh Noughabi, H. (2010). A new estimator of entropy and its application in testing normality.\emph{ Journal of Statistical Computation and Simulation}, 80, 1151--1162.

\bibitem {bib13} Alizadeh Noughabi, H., and Arghami, N. R.  (2010). A new estimator of entropy. \emph{Journal of the Iranian Statistical Society}, 9, 53--64.

\bibitem {bib4} Al-Labadi, L. (2018). On Metrizing Vague Convergence of Random Measures with Applications on Bayesian Nonparametric Models. \emph{Statistics}, 52, 445--457.

\bibitem {bib5}  Al-Labadi, L., and Abdelrazeq, I.  (2017). On Functional Central Limit Theorems of  Bayesian Nonparametric Priors.  \emph{Statistical Methods \& Applications}, 26,  215--229.

\bibitem {bib6} Al-Labadi, L., and Evans, M. (2018). Prior-Based Model Checking. \emph{Canadian Journal of Statistics}, 46, 380--398.

\bibitem {bib7}Al-Labadi, L., and Zarepour, M. (2014a). Goodness of fit tests
based on the distance between the Dirichlet process and its base measure.
\emph{Journal of Nonparametric Statistics}, 26, 341--357.

\bibitem {bib8} Al-Labadi, L., and Zarepour, M. (2014b). On Simulations from the Two-Parameter Poisson-Dirichlet Process and the Normalized Inverse-Gaussian Process. \emph{Sankhy\=a  A}, 76,  158--176.

\bibitem {bib9} Al-Labadi, L., and Zarepour, M. (2013a).  A Bayesian Nonparametric Goodness of Fit Test for Right Censored Data Based on Approximate Samples from the Beta-Stacy Process. \emph{Canadian Journal of Statistics}, 41,  466--487.

\bibitem {bib10} Al-Labadi, L., and Zarepour, M. (2013b). On Asymptotic Properties and Almost Sure Approximation of the Normalized Inverse-Gaussian Process.  \emph{Bayesian Analysis}, 8,  553--568.

\bibitem {bib11} Al-Omari, A. I. (2014). Estimation of entropy using random sampling. \emph{Journal of Computation and Applied Mathematics}, 261, 95--102.

\bibitem {bib12} Al-Omari, A. I. (2016).  A new measure of entropy of continuous random variable. Journal of Statistical Theory and Practice, 10, 721--735

\bibitem {bib28}  Atkinson, C. M.  (2017). Weak and strong laws of large numbers of negatively associated random variables. PhD thesis, University of Regina. \url{https://ourspace.uregina.ca/bitstream/handle/10294/7844/Atkinson_Christopher_200324544_PHD_MATH_Fall2017.pdf}

\bibitem {bib13} Bhattacharya, G., Ghosh, K.,  and Chowdhur, A. S. (2012). An affinity-based new local distance function and similarity measure for kNN algorithm. \emph{Pattern Recognition Letters}, 33, 356--363.

\bibitem {bib14} Beirlant, J., Dudewicz, E. J., Gy\"oria , L., and van der Meulen, E. C. (1997). Nonparametric entropy estimation: an overview. \emph{International Journal of Mathematics and Statistics}, 6,  17--39.


\bibitem {bib15} Berrett, T. B.,   Grose, D. J., and Samworth, R. J. (2018). Package `IndepTest’: Nonparametric Independence Tests Based on Entropy Estimation. https://cran.r-project.org/web/packages/IndepTest/IndepTest.pdf


\bibitem {bib16} Berrett, T. B., Samworth, R. J., and Yuan, M. (2019). Efficient multivariate entropy estimation via
k-nearest neighbour distances. \emph{Annals of Statistics}, 47, 288--318.




\bibitem {bib17} Bondesson, L. (1982). On simulation from infinitely divisible
distributions. \emph{Advances in Applied Probability}, {14}, 885--869.

\bibitem {bib18} Bouzebda, S., Elhattab, I.,  Keziou, A., and Lounis, T.  (2013). New Entropy Estimator with an Application to Test of Normality.\emph{ Communications in Statistics - Theory and Methods}, 42, 2245--2270.



\bibitem {bib19}Correa, J. C. (1995).  A new estimator of entropy. \emph{Communications
in Statistics - Theory and Methods}, 24, 2439--2449.

\bibitem {bib20}Cover, T. M., and  Thomas, J. A. (1991).  \emph{Elements of Information Theory}, Wiley.


\bibitem {bib21} Ebrahimi, N., Pflughoeft, K., and Soofi, E. (1994). Two measures of sample entropy. \emph{Statistics \& Probability Letters}, 20, 225--234.


\bibitem {bib22}Ferguson, T. S. (1973). A Bayesian analysis of some
nonparametric problems. \emph{Annals of Statistics}, {1}, 209--230.

\bibitem {bib23}Ferguson, T. S. (1996). \emph{A Course in Large Sample Theory}.  Chapman \& Hall/CRC.


\bibitem {bib24} Grzegorzewski, P. and Wieczorkowski, R. (1999). Entropy-based goodness-of-fittest for exponentiality. \emph{CommunicationsinStatistics - Theory and Methods}, 28, 1183--1202.


\bibitem {bib25} Ishwaran, H., and James, L. F. (2001). Gibbs Sampling Methods for Stick-Breaking
Priors. \emph{Journal of the American Statistical Association}, 96, 161--173.

\bibitem {bib26}Ishwaran, H., and Zarepour, M. (2002). Exact and Approximate
Sum Representations for the Dirichlet Process. \emph{The Canadian Journal of
Statistics}, 30, 269--283.




\bibitem {bib27} Kozachenko, L. F., and Leonenko, N. N. (1987). Sample estimate of the entropy of a random
vector. \emph{Problems of Information Transmission},  23,  95--101.


\bibitem {bib28}Mazzuchi, T. A., Soofi, E. S., and Soyer, R. (2000). Computations of maximum entropy Dirichlet for
modeling lifetime data. \emph{Computational Statistics and Data Analysis}, 32, 361--378.

\bibitem {bib29}Mazzuchi, T. A., Soofi, E. S., Soyer, R. (2008). Bayes estimate and inference for entropy
and information index of fit. \emph{Econometric Reviews}, 27(4--6), 428--456.

\bibitem {bib30} Mitra, P., Murthy, C. A., and  Pal, S. K. (2002). Unsupervised feature selection using feature similarity. \emph{IEEE Transactions on Pattern Analysis and Machine Intelligence}, 24, 301--312.



\bibitem {bib31}Sethuraman, J. (1994). A constructive definition of Dirichlet
priors. \emph{Statistica Sinica}, {4}, 639--650.

\bibitem {bib32}Shannon, C. E. (1948). A mathematical theory of communication. The Bell System Technical Journal, 27, 379--423, 623--656.

\bibitem {bib33} Vasicek, O. (1976). A test for normality based on sample entropy. \emph{Journal of Royal Statistical Society B}, 38, 54--59.

\bibitem {bib34}van Es, B. (1992). Estimating functionals related to a density
by a class of statistics based on spacings. \emph{Scandinavian Journal of Statistics}, 19,  61--72.


\bibitem {bib35}Wolpert, R. L., and Ickstadt, K., (1998). Simulation of
L\'{e}vy random fields. In \emph{Practical Nonparametric and Semiparametric
Bayesian Statistics}, ed. D. Day, P. Muller, and D. Sinha, Springer, 227--242.

\bibitem {bib36} Wieczorkowski, R. and Grzegorzewski, P. (1999). Entropy estimators-improvements and comparisons. \emph{Communications in Statistics - Simulation and Computation}, 28, 541--567.


\bibitem {bib37}Zarepour, M., and Al-Labadi, L. (2012). On a rapid simulation
of the Dirichlet process. \emph{Statistics \& Probability Letters}, 82,  916--924.
\end{thebibliography}
\end{document}